\documentclass[a4paper]{amsart}
\usepackage[all]{xy}
\usepackage{amssymb}
\usepackage{amsmath}
\usepackage{bm}
\usepackage[colorlinks,linkcolor=blue]{hyperref}
\usepackage{color,xcolor,dsfont}
\usepackage{enumerate}
\usepackage{epsfig}
\usepackage{extarrows}
\usepackage{framed}
\usepackage{graphicx}
\usepackage{marvosym}
\usepackage{mathrsfs}
\usepackage{setspace} 
\usepackage{tikz}
\usepackage{tikz-cd}

\newcommand{\scha}{\mathcal{A}}

\newcommand{\schf}{\mathcal{F}}
\newcommand{\schp}{\mathcal{P}}
\newcommand{\schs}{\mathcal{S}}
\newcommand{\scht}{\mathcal{T}}

\newcommand{\she}{\mathscr{E}}
\newcommand{\shf}{\mathscr{F}}
\newcommand{\shg}{\mathscr{G}}
\newcommand{\shh}{\mathscr{H}}
\newcommand{\shi}{\mathscr{I}}

\newcommand{\shl}{\mathscr{L}}
\newcommand{\shm}{\mathscr{M}}

\newcommand{\sho}{\mathscr{O}}

\newcommand{\sht}{\mathscr{T}}

\DeclareMathOperator{\Pic}{Pic}
\DeclareMathOperator{\Aut}{Aut}

\DeclareMathOperator{\Hom}{Hom}

\DeclareMathOperator{\rank}{rank}

\DeclareMathOperator{\ch}{ch}
\DeclareMathOperator{\CH}{CH}

\DeclareMathOperator{\Geo}{Geo}
\DeclareMathOperator{\Stab}{Stab}
\DeclareMathOperator{\Slic}{Slic}

\DeclareMathOperator{\gldim}{gldim}

\newcommand{\bz}{\mathbb{Z}}
\newcommand{\bh}{\mathbb{H}}
\newcommand{\bq}{\mathbb{Q}}
\newcommand{\br}{\mathbb{R}}
\newcommand{\bc}{\mathbb{C}}

\newcommand{\bp}{\mathbb{P}}

\newcommand{\DCoh}{\cate{D}^b}

\DeclareMathOperator{\GL}{GL}
\newcommand{\GLp}{\operatorname{GL^+}}
\newcommand{\grp}{{\tilde{\GLp}}(2,\br)}

\newcommand{\defi}[1]{{\textbf{\emph{#1}}}}
\newcommand{\cate}[1]{{\textbf{\textup{#1}}}}

\theoremstyle{plain}
\newtheorem{theorem}{Theorem}[section]
\newtheorem{corollary}[theorem]{Corollary}
\newtheorem{lemma}[theorem]{Lemma}
\newtheorem{proposition}[theorem]{Proposition}

\theoremstyle{definition}
\newtheorem{definition}[theorem]{Definition}
\newtheorem{remark}[theorem]{Remark}
\newtheorem{example}[theorem]{Example}

\author{Ziqi Liu}
\title{Stability conditions on irreducible projective curves}
\address{Dipartimento di Matematica “F. Enriques”, Università degli Studi di Milano, Via Cesare Saldini 50, 20133 Milano, Italy.}
\email{ziqi.liu@unimi.it}



\thanks{The author is a member of the INdAM group GNSAGA (2025).}

\begin{document}

\maketitle
\thispagestyle{plain}

\begin{abstract}
This note revisits stability conditions on the bounded derived categories of coherent sheaves on irreducible projective curves. In particular, all stability conditions on smooth curves are classified and a connected component of the stability manifold containing all the geometric stability conditions is identified for singular curves. On smooth curves of positive genus, the set of all non-locally-finite stability conditions gives a partial boundary of any known compactification of the stability manifold. To provide a reasonable full boundary, a notion of regular weak stability condition is proposed based on the definition of Collins--Lo--Shi--Yau and is classified for smooth curves of positive genus. On non-rational singular curves, any locally-finite numerical stability condition is shown to be geometric.
\end{abstract}

\section{Introduction}
The notion of stability condition on a triangulated category is introduced in the seminal article \cite{Bri07} by Bridgeland as a mathematical formulation for Douglas' ideas of stability in string theory. 

A stability condition $\sigma=(\schp,Z)$ on a triangulated category $\scht$ consists of a collection $\schp$ of additive subcategories $\schp(\phi)\subset\scht$ for any real number $\phi\in\br$ and a group homomorphism $Z$ from the Grothendieck group $K(\scht)$ of $\scht$ to $\bc$ such that some compatible conditions are satisfied. 

Given any interval $I\subset\br$, the subcategory $\schp(I)$ generated by $\schp(\phi)$ for all $\phi\in I$ is quasi-abelian. A stability condition $(\schp,Z)$ is called \emph{locally finite} if there exists $\eta>0$ such that the quasi-abelian category $\schp((\phi-\eta,\phi+\eta))$ is of finite length for any $\phi$. The set of all locally finite stability conditions on a triangulated category $\scht$ admits a canonical complex structure, once it is non-empty, and the resulting complex manifold is denoted by $\Stab(\scht)$. 

This note studies the complex manifold $\Stab(\DCoh(C))$ associated with the bounded derived category $\DCoh(C)$ of coherent sheaves on an irreducible projective curve $C$.

\subsection{Stability conditions on smooth curves} Let $C$ be a smooth projective curve, then one defines the numerical Grothendieck group $N(C):=K(C)/\ker(\chi)$ by quotient out numerically trivial classes. The complex manifold $\Stab(\DCoh(C))$ has a finite dimensional submanifold $\Stab(C)$ containing all the locally finite numerical stability conditions on $C$. Here a stability condition $\sigma=(\schp,Z)$ on $\cate{D}^b(C)$ is called \emph{numerical} when the central charge $Z\colon K(C)\rightarrow\bc$ factors through the numerical Grothendieck group $N(C)$ via the quotient map $K(C)\twoheadrightarrow N(C)$. 

The complex submanifold $\Stab(C)$ has been identified in \cite{Oka06} for $C\cong\bp^1$, in \cite{Bri07} for elliptic curves, and in \cite{Mac07} for smooth curves of genus $g\geq 2$. However, the whole space $\Stab(\DCoh(C))$ is only known when $C\cong\bp^1$ as $K(\bp^1)=N(\bp^1)$. It turns out that it is in fact true for all smooth projective curves.

\begin{theorem}\label{1_main-1-1}
Consider a smooth curve $C$ of positive genus, then $\Stab(\DCoh(C))=\Stab(C)\cong\bc\times\bh$ as a complex manifold where $\bh$ is the hyperbolic upper half-plane.
\end{theorem}

There are no non-locally-finite stability conditions on $\bp^1$ by \cite{GKR04,Oka06}. This note classifies non-locally-finite stability conditions on other smooth curves.

\begin{theorem}\label{1_main-1-2}
The set of non-locally-finite stability conditions on $\DCoh(C)$ is equal to $\bc\times(\br-\bq)$ for a smooth curve $C$ of positive genus.
\end{theorem}

In particular, up to the action of a complex number, a non-locally-finite stability condition on $\DCoh(C)$ can be identified with a pair $\sigma_{\beta}=(\schp_{\beta},Z_{\beta})$ for some non-rational real number $\beta\in\br-\bq$. The same construction instead gives a locally finite stability condition $\sigma_{\beta}$ on $\DCoh(\bp^1)$ for any $\beta\in\br-\bz$.

\subsection{Boundary points of the complex manifold}
As pointed out by Bridgeland in his ICM talk \cite[Remark 4.2]{Bri06}, one should expect that the non-locally-finite stability conditions on $\scht$ are boundary points for $\Stab(\scht)$. There are different approaches on the compactification of $\Stab(\scht)$, such as \cite{BDL23,BMS24,Bol23,CLSY25,De25,KKO24}. 

This expectation will be confirmed for the complex manifold $\Stab(C)\cong\bc\times\bh$ for smooth curves with positive genus. It turns out that each $\sigma_{\beta}$ is an element in the horizontal real axis with respect to $\bh$ for the central charge partial compactification \cite{Bol23} of $\Stab(C)$ and the Thurston compactification of the quotient $\Stab(C)/\bc\cong\bh$ in \cite{BDL23,KKO24}. Moreover, they can also be seen as the boundary points with respect to the global dimension defined in \cite{IQ23,Qiu23}. 

A rational number $\beta$ corresponds to infinite many pairs $(\schp,Z_{\beta})$ which fail to meet the definition of a stability condition. So to fill the missing points in an appropriate way, one needs a suitable weaker notion of stability. In this note, we suggest a definition of weak stability conditions based on \cite{CLSY25}. This notion provides, in many senses, the whole boundary for $\Stab(C)$ or $\Stab(C)/\bc$ for a positive genus curve $C$ such that a $\beta$ only corresponds to two pairs $(\schp_{\beta},Z_{\beta})$ and $(\schp'_{\beta},Z_{\beta})$. The weak stability conditions are also classified for $\DCoh(C)$.

\subsection{Stability conditions on singular curves}
In the end, we will investigate stability conditions on singular curves. The first obstruction is that the Grothendieck group of $\DCoh(C)$ is unclear in this case. So one starts with a closed submanifold $\Geo^{\dag}(C)\cong\bc\times\bh$ in $\Stab(\DCoh(C))$ containing stability conditions induced by slope stability of coherent sheaves. It has been proved in \cite{BK06} for singular Weierstraß cubic curves that $\Geo^{\dag}(C)=\Stab(\DCoh(C))$. In general, one can show 

\begin{theorem}\label{1_main-3}
Consider a singular curve $C$, then the closed submanifold $\Geo^{\dag}(C)$ is a connected component of $\Stab(\DCoh(C))$ and contains all the stability conditions on $\DCoh(C)$ such that the skyscraper sheaves are stable with the same phase (i.e.\ the geometric stability conditions).
\end{theorem}

Inspired by \cite{FLZ22}, we are able to show that every locally finite numerical stability condition is geometric as $\Pic^0(C)$ is always a smooth algebraic group.

\begin{proposition}
Consider a non-rational curve $C$, then $\Stab(C)=\Geo^{\dag}(C)$.
\end{proposition}

It is expected that one has $\Stab(\DCoh(C))=\Geo^{\dag}(C)$. However, it requires a better understanding of the possible bounded $t$-structures on $\DCoh(C)$.

\subsection*{Conventions}
In this note, everything is over the field of complex numbers, curves are integral and projective, categories are essentially small, and functors are derived.

\section{Backgrounds on Stability Conditions}
\subsection{An elementary fact about $t$-structures}
At first, we recall an elementary fact about bounded $t$-structures which lacks clear references. The definitions and terminology are cited or borrowed from \cite{GM}.

\begin{definition}
A \defi{bounded $t$-structure} on a triangulated category $\scht$ consists of a pair of strictly full subcategories $(\scht^{\leq0},\scht^{\geq0})$ such that
\begin{itemize}
	\item $\scht=\bigcup_{m,n\in\bz}\scht^{\geq m}\cap\scht^{\leq n}$;
	\item $\scht^{\leq0}\subseteq\scht^{\leq1}$ and $\scht^{\geq1}\subseteq\scht^{\geq0}$;
	\item $\Hom(X,Y)=0$ for $X\in\scht^{\leq0}$ and $Y\in\scht^{\geq 1}$;
	\item for any non-zero $X\in\scht$ there is an exact triangle $$A\rightarrow X\rightarrow B\rightarrow A[1]$$ 
	with $A\in\scht^{\leq0}$ and $B\in\scht^{\geq1}$;
\end{itemize}
where $\scht^{\leq n}:=\scht^{\leq0}[-n]$ and $\scht^{\geq n}:=\scht^{\geq0}[-n]$ for any $n\in\bz$. 
\end{definition}

The full subcategory $\scha=\scht^{\geq0}\cap\scht^{\leq0}\subseteq\scht$ is called the \emph{heart} of the bounded $t$-structure and it is always an abelian category \cite[IV.4.4]{GM}. Moreover, one has

\begin{lemma}[{{\cite[Lemma 3.2]{Bri07}}}]\label{2_lemma-t-structure}
A full additive subcategory $\scha$ of a triangulated category $\scht$ is the heart of a bounded $t$-structure on $\scht$ if and only if 1) for every non-zero object $E\in\scht$ there exists a finite sequence of integers $k_1>k_2>\cdots>k_n$ and a collection of triangles
\begin{displaymath}
	0= \xymatrix@R=2em{
		E_0\ar[r]&E_1\ar[r]\ar[d]&E_2\ar[r]\ar[d]&\cdots\ar[r]&E_{n-1}\ar[r]\ar[d]&E_n\ar[d]\\
		&A_1\ar@{-->}[ul]&A_2\ar@{-->}[ul]&\ar@{-->}[ul]\cdots&\cdots&A_n\ar@{-->}[ul]
	} =E
\end{displaymath}
with $A_i\in\scha[k_i]$ for all $i$; and 2) one has $$\Hom(A[m_1],B[m_2])=0$$ 
for any $A,B\in\scha$ and $m_1>m_2$.
\end{lemma}

In fact, the sequence $k_1>\cdots>k_n$ is unique for a non-zero object $E\in\scht$ once a bounded $t$-structure $(\scht^{\geq0},\scht^{\leq0})$ is fixed. Also, one has $E\in\scht^{\geq-k_1}-\scht^{\geq-k_1+1}$ and $E\in\scht^{\leq-k_n}-\scht^{\leq-k_n-1}$.

\begin{lemma}\label{2_fact-t-structure}
Consider two bounded $t$-structures $(\scht_1^{\leq0},\scht_1^{\geq0})$ and $(\scht_2^{\leq0},\scht_2^{\geq0})$ on a triangulated category $\scht$ with hearts $\scha_1,\scha_2$ respectively, then $(\scht_1^{\leq0},\scht_1^{\geq0})$ is equal to $(\scht_2^{\leq0},\scht_2^{\geq0})$ once $\scha_1\subseteq\scha_2$.

\begin{proof}
At first, one observes that $\scha_1=\scha_2$ otherwise for any $A\in\scha_2-\scha_1$ one can take integers $k_1>\cdots >k_n$ and two non-trivial morphisms $A_1\rightarrow A$ and $A\rightarrow A_n$ with $A_1\in\scha_1[k_1]\subseteq\scha_2[k_1]$ and $A_n\in\scha_1[k_n]\subseteq\scha_2[k_n]$ by Lemma \ref{2_lemma-t-structure}. It follows an inequality $0\geq k_1>k_n\geq 0$, absurd! 

Now we are ready to prove $(\scht_1^{\leq0},\scht_1^{\geq0})=(\scht_2^{\leq0},\scht_2^{\geq0})$. Since $\scha_1=\scha_2$ one can obtain $\scht_1^{\leq0}=\scht_2^{\leq0}$ and $\scht_1^{\geq0}=\scht_2^{\geq0}$ by Lemma \ref{2_lemma-t-structure} without much effort.
\end{proof}
\end{lemma}

\subsection{Stability condition on triangulated categories}
Here we introduce basic notions about stability conditions on triangulated categories following \cite{Bri07}.

\begin{definition}
A \defi{slicing} $\schp$ on a triangulated category $\scht$ is a collection of full additive subcategories $\schp(\phi)\subset\scht$ for each $\phi\in\br$, satisfying the following axioms:
\begin{itemize}
	\item for all $\phi\in\br$, one has $\schp(\phi+1)=\schp(\phi)[1]$;
	\item one has $\Hom(\schp(\phi_1),\schp(\phi_2))=0$ whenever $\phi_1>\phi_2$;
	\item for any $0\neq E\in\scht$, one has a finite sequence of real numbers $\phi_1>\phi_2>\cdots>\phi_n$ and a collection of triangles
	\begin{displaymath}
		0= \xymatrix@R=2em{
			E_0\ar[r]&E_1\ar[r]\ar[d]&E_2\ar[r]\ar[d]&\cdots\ar[r]&E_{n-1}\ar[r]\ar[d]&E_n\ar[d]\\
			&A_1\ar@{-->}[ul]&A_2\ar@{-->}[ul]&\ar@{-->}[ul]\cdots&\cdots&A_n\ar@{-->}[ul]
		} =E
	\end{displaymath}
with $A_i\in\schp(\phi_i)$ for all $i$. 
\end{itemize}
The filtration is called a \emph{Harder--Narasimhan filtration} of $E$ with respect to $\schp$.
\end{definition}

The subcategory $\schp(\phi)$ is always quasi-abelian. Moreover, one defines $\schp(I)$, for an interval $I\subset\br$, to be the extension-closed subcategory of the triangulated category $\scht$ generated by the subcategories $\schp(\phi)$ for $\phi\in I$. Then $\schp((0,1])$ is an abelian category and is the heart of a bounded $t$-structure on $\scht$. So the category $\schp((0,1])$ is called the \emph{heart} of the slicing $\schp$.

\begin{definition}
A \defi{pre-stability condition} $\sigma=(\schp,Z)$ on a triangulated category $\scht$ consists of a group homomorphism $Z\colon K(\scht)\rightarrow\bc$ (called the \emph{central charge} of the pre-stability condition $\sigma$) and a slicing $\schp$ such that $Z(E)=m(E)\exp(\bm{i}\pi\phi)$ with $m(E)\in\br_{\geq 0}$ for any $0\neq E\in\schp(\phi)$. 
\end{definition}

\begin{remark}
One notices that this definition is different from the notion of weak pre-stability condition in \cite{PT19} as the phase of an object in $\schp((0,1])$ with zero central charge might not be $1$.
\end{remark}

A non-zero object in $\schp(\phi)$ is called \emph{semistable of phase $\phi$} with respect to $\sigma$ and the simple objects of $\schp(\phi)$ are called \emph{stable of phase $\phi$} with respect to $\sigma$. Similar to the case of bounded $t$-structures, the decomposition for a non-zero object $E$ is unique up to isomorphism and the objects $A_i$ are called the \emph{semistable factors} of the object $E$ with respect to $\sigma$. Also, one could define the real numbers $\phi_+(E)=\phi_1$ and $\phi_-(E)=\phi_n$. By definition, one sees

\begin{proposition}
Consider a pre-stability condition $\sigma=(\schp,Z)$ on $\scht$, then for any $E\neq0$ in $\schp((0,1])$, one has $Z(E)\in\mathbb{U}\cup\br_{\leq 0}$ where $\mathbb{U}$ is the upper-half plane.
\end{proposition}

Similar to \cite[Lemma 8.2]{Bri07}, there are two mutually-commutative group actions on the set of all pre-stability conditions on a given triangulated category $\scht$.

The group $\Aut(\scht)$ of exact autoequivalences on $\scht$ acts on the left, via
$$\Phi.(\schp,Z)=(\schp',Z\circ\Phi_*^{-1})$$
where $\Phi_*$ denotes the induced automorphism of $K(\scht)$ and $\schp'(\phi):=\Phi(\schp(\phi))$.

The universal cover 
$$
\grp  = \left \{ \tilde{g}=(M, f) ~  \middle\vert ~ 
\begin{array}{l} 
	M \in \GLp(2,\br), f \colon \br \rightarrow \br \text{ is an increasing function} \\ 
	\text{such that for all $\phi \in \br$ we have $f(\phi+1) = f(\phi) + 1$} \\
	\text{and $M \cdot e^{\bm{i}\pi \phi} \in \br_{>0} \cdot e^{\bm{i} \pi f(\phi)}$} 
\end{array}  \right \}
$$
of the group $\GLp(2,\br)=\{M\in\GL(2,\br)\,|\,\det(M)>0\}$ acts on the right via 
$$
(\schp, Z) . (M, f):= (\schp', M^{-1} \circ Z). 
$$ 
where $\schp'(\phi):=\schp(f(\phi))$ and $\bc=\br\oplus\br \bm{i}$ is used to validate $M^{-1}\circ Z\colon K(\scht)\rightarrow\bc$.

There exists a natural subgroup $\bc\subset\grp$ such that for any $\lambda\in\bc$ one has
$$(\schp,Z).\lambda=(\schp',\exp(-\bm{i}\pi\lambda)\cdot Z)$$
where $\schp'(\phi):=\schp(\phi+\Re(\lambda))$ and the multiplication with $Z$ is the usual one. 

\begin{definition}
A \defi{stability condition} $\sigma=(\schp,Z)$ on a triangulated category $\scht$ is a pre-stability condition on $\scht$ such that $Z(E)\neq0$ for any $0\neq E\in\schp(\phi)$.
\end{definition}

The two group actions about pre-stability conditions restricts to the set of all the stability conditions on $\scht$. Moreover, one has $Z(E)\in\mathbb{U}\cup\br_{< 0}$ for a stability condition $\sigma=(\schp,Z)$ and a non-zero object $ E\in\schp((0,1]))$.

\subsection{The space of locally finite stability conditions}
One of the main results in \cite{Bri07} is to construct a canonical complex structure on the set $\Stab(\scht)$ of locally finite stability conditions. 

\begin{definition}
A pre-stability condition $\sigma=(\schp,Z)$ is called \defi{locally finite} if there exists some $\epsilon>0$ such that the (quasi-abelian) category $\schp((\phi-\varepsilon,\phi+\varepsilon))$ is of finite length for any $\phi\in\br$.
\end{definition}

The quasi-abelian category $\schp(\phi)$ associated with a locally finite pre-stability condition $\sigma=(\schp,Z)$ is of finite length as well so that every object in $\schp(\phi)$ has a finite Jordan--Hölder filtration into stable factors (which are precisely the simple objects by definition) of the same phase.

\begin{theorem}{{\cite[Theorem 1.2]{Bri07}}}\label{2_main-Bridgeland}
Consider a triangulated category $\scht$ and a connected component $\Sigma\subset\Stab(\scht)$, then there exists a linear subspace $V(\Sigma)$ of $\Hom_{\bz}(K(\scht),\bc)$ with a well-defined linear topology such that the natural map $\mathcal{Z}\colon\Sigma\rightarrow V(\Sigma)$ defined by $(\schp,Z)\mapsto Z$ is a local homomorphism. 
\end{theorem}

\begin{remark}
The two group actions described in the previous subsection restricts to group actions on the generalized metric space $\Stab(\scht)$.
\end{remark}

There are two types of special pre-stability conditions on the bounded derived category $\DCoh(X)$ of a smooth variety $X$.

\begin{definition}
A pre-stability condition $\sigma$ on $\DCoh(X)$ is called \defi{numerical} if its central charge $Z$ factors through the numerical Grothendieck group $N(X)$ via the quotient map $K(X)\twoheadrightarrow N(X)$.
\end{definition}

The two group actions on the set of all pre-stability conditions preserve the subset of all the numerical ones. A numerical pre-stability condition is not necessary locally finite and the subspace $\Stab(X)\subset\Stab(\DCoh(X))$ of numerical locally finite stability conditions is a finite dimensional complex manifold with the subspace topology.

\begin{definition}
A pre-stability condition $\sigma$ on $\DCoh(X)$ is called \defi{geometric} if all the skyscraper sheaves on $X$ are $\sigma$-semistable with the same phase.
\end{definition}

The geometric pre-stability conditions are only preserved by the group action of $\grp$ and is not a priori numerical or locally finite. One uses $\Geo(X)$ to denote the subspace of all geometric locally finite stability conditions on $\DCoh(X)$.

\subsection{Some pre-stability conditions on curves}
Let $C$ be a curve, and let $K(C):=K(\DCoh(C))$ be its Grothendieck group. Then the rank and degree (defined as $\chi(\shf)$) of coherent sheaves induce a surjection $(\rank,\deg)\colon K(C)\rightarrow\bz\oplus\bz$ which is an isomorphism for rational curves. 

\begin{definition}
A stability condition $\sigma$ on $\DCoh(C)$ is called \defi{numerical} if its central charge factors through the surjection $(\rank,\deg)\colon K(C)\rightarrow\bz^{\oplus 2}$.
\end{definition}

The following (pre)-stability conditions on $\DCoh(C)$ are based on the existence of Harder--Narasimhan filtrations and Jordan--Hölder filtrations for the slope stability of coherent sheaves (see \cite{HL} for more details).

\begin{example}[{{\cite[Example 5.4]{Bri07} and \cite{BK06}}}]\label{2_slope-stability}
The slope stability gives a locally finite numerical stability condition $(\schp,Z)$ on $\DCoh(C)$ such that $\schp(\phi)$ contains the zero object and all slope semistable coherent sheaves on $C$ with phase $\phi$ with respect to the slope stability for $\phi\in(0,1]$, and the central charge $Z(r,d)=-d+r\bm{i}$.
\end{example}

\begin{example}\label{2_slope-stability-degenerate}
Let $\beta\in\br\cup\{\infty\}$, then one has a torsion pair 
\begin{align*}
		\schf_{\beta}&=\{\shg\in\cate{Coh}(C)\,|\,\shg\textup{ is torsion-free and }\mu_+(\shg)< \beta\}\\
		\scht_{\beta}&=\{\shg\in\cate{Coh}(C)\,|\,\shg\textup{ is torsion or its torsion-free part satisfies }\mu_-(\shg)\geq \beta\}
\end{align*}
on $\cate{Coh}(C)$, where $\mu_+(\shg)$ (resp.\ $\mu_-(\shg)$) denotes the maximal (resp.\ minimal) slope of the semistable factors of $\shg$. It follows a bounded $t$-structure on $\DCoh(C)$ with heart
$$\scha_{\beta}=\{A\in\DCoh(C)\,|\,\shh^{-1}(A)\in\schf_{\beta},\shh^0(A)\in\scht_{\beta},\shh^i(A)=0,i\neq0,-1\}$$
according to \cite{HRS}. One notices that the heart $\scha_{\beta}$ is denoted by $\textsf{\textup{A}}(\theta,\textsf{\textup{P}}(\theta)^s)$ in the statement \cite[Proposition 5.1]{BK06} where $\theta\in(0,1]$ is the usual phase associated with slope stability for slope $\beta$. Then one has a numerical pre-stability condition $\sigma_{\beta}=(\schp_{\beta},Z_{\beta})$ on $\DCoh(C)$ such that $\schp_{\beta}(1)=\schp_{\beta}((0,1])=\scha_{\beta}$ and the central charge of it is $Z_{\beta}(r,d)=-d+\beta r$ for $\beta\in\br$ or $Z_{\infty}(r,d)=r$. 
\end{example}

One notices that $\schf_{\infty}$ contains all the torsion-free sheaves on $C$ and $\scht_{\infty}$ contains all the torsion sheaves on $C$. 

\begin{example}\label{2_slope-stability-degenerate-alternative}
Let $\beta\in\br$, then one has another torsion pair 
\begin{align*}
		\schf'_{\beta}&=\{\shg\in\cate{Coh}(C)\,|\,\shg\textup{ is torsion-free and }\mu_+(\shg)\leq \beta\}\\
		\scht'_{\beta}&=\{\shg\in\cate{Coh}(C)\,|\,\shg\textup{ is torsion or its torsion-free part satisfies }\mu_-(\shg)>\beta\}
\end{align*}
on $\cate{Coh}(C)$. It follows a bounded $t$-structure on $\DCoh(C)$ with heart
$$\scha'_{\beta}=\{A\in\DCoh(C)\,|\,\shh^{-1}(A)\in\schf'_{\beta},\shh^0(A)\in\scht'_{\beta},\shh^i(A)=0,i\neq0,-1\}$$
according to \cite{HRS}. One notices that the heart $\scha_{\beta}'$ is denoted by $\textsf{\textup{A}}(\theta,\varnothing)$ in the statement \cite[Proposition 5.1]{BK06} where $\theta\in(0,1]$ is the usual phase associated with slope stability for slope $\beta$. Also, the heart $\scha'_{\beta}$ is equal to $\scha_{\beta}$ once there are no semistable coherent sheaves on $C$ of slope $\beta$. Then one has a numerical pre-stability condition $\sigma'_{\beta}=(\schp'_{\beta},Z_{\beta})$ on $\DCoh(C)$ with $\schp'_{\beta}(1)=\schp'_{\beta}((0,1])=\scha'_{\beta}$. 

Moreover, one defines a pre-stability condition $\sigma'_{\infty}=(\schp'_{\infty},Z_{\infty})$ where the slicing $\schp'_{\infty}$ is determined by $\schp'_{\infty}(0)=\cate{Coh}(C)$. So one has $\sigma'_{\beta}$ for any $\beta\in\br\cup\{\infty\}$.
\end{example}

One notices that $\sigma_{\beta}=\sigma'_{\beta}$ when there are no semistable sheaves on $C$ with slope $\beta$. In this case, it is a stability condition on $\DCoh(C)$ similar to \cite[Example 5.6]{Bri07}.

\section{Stability conditions on smooth curves of positive genus}
In this section, we will revisit the stability conditions on the bounded derived category of a smooth curve of positive genus.  

\subsection{The complex manifold of numerical stability conditions}
Let us fix a positive genus smooth curve $C$ in this section. At first, one recalls:

\begin{lemma}[{{\cite[Lemma 7.2]{GKR04}}}]\label{3-key-technical-lemma}
Consider a coherent sheaf $\shg$ on $C$, then for any exact triangle 
$$X\rightarrow \shg\rightarrow Y\rightarrow X[1]$$ 
in $\DCoh(C)$ with $\Hom(X,Y[n])=0$ for $n\leq0$ the objects $X,Y$ are coherent sheaves.
\end{lemma}

This lemma does not hold for $\bp^1$ and is not known for singular curves. It follows a critical corollary by the proof for \cite[Theorem 2.7]{Mac07}.

\begin{corollary}\label{3-lemma-macri}
The skyscraper sheaves and line bundles are semistable with respect to any pre-stability condition $\sigma$ on $\DCoh(C)$. They are stable if $\sigma$ is also locally finite. Moreover, one has
$$\phi(\sho_x)-1\leq\phi(\shg)\leq\phi(\sho_x)$$ 
for any $\sigma$-semistable locally free sheaf $\shg$ on $C$ and any $x\in C$.
\end{corollary}

Therefore, a locally finite numerical pre-stability condition on $\DCoh(C)$ is geometric and equal to the one in Example \ref{2_slope-stability} up to a unique $\grp$-action.

\begin{theorem}[{{\cite[Theorem 2.7]{Mac07}}}]\label{3-macri-main-result}
The action of $\grp$ on $\Stab(C)$ is free and transitive.
\end{theorem}

It is proved in \cite{Woolf12} that as a complex manifold one has $\Stab(C)\cong\bc\times\bh$ where $\bh$ is the hyperbolic upper half-plane. In particular, a point $(\lambda,\beta+\alpha\bm{i})\in\bc\times\bh$ corresponds to the stability condition $\sigma_{\alpha,\beta}.\lambda$ according to \cite[Theorem 2.7]{Mac07} and the arguments in \cite[Section 3]{BMW15}, where $\sigma_{\alpha,\beta}$ is the stability condition with central charge $Z_{\alpha,\beta}(r,d)=-d+(\beta+\alpha\bm{i})r$ and heart $\cate{Coh}(C)$.

\begin{remark}
Since $\Pic(C)$ is generated by the classes of points on $C$, a geometric pre-stability condition on $\DCoh(C)$ is numerical. So one has $\Geo(C)=\Stab(C)$.
\end{remark}

\subsection{Non-existence of non-numerical stability conditions}
Now we are ready to prove the first main result of this note. One recalls that the Grothendieck group $K(C)$ is isomorphic to the direct sum $\bz\oplus\Pic(C)$ and the Picard group $\Pic(C)$ is non-canonically isomorphic to $\bz\oplus\Pic^{\circ}(C)$ such that the subgroup $\Pic^{\circ}(C)$ contains the classes of degree zero divisors.

\begin{proposition}\label{3-main-1-1}
There are no non-numerical pre-stability conditions on $\DCoh(C)$.
\end{proposition}

\begin{proof}
Otherwise, one chooses a non-numerical pre-stability condition $\sigma=(\schp,Z)$ on $\DCoh(C)$. At first, due to Corollary \ref{3-lemma-macri}, the supremum $\phi_1=\sup\{\phi(\sho_x)\,|\,x\in C\}$ and the infimum $\phi_2=\inf\{\phi(\sho_x)\,|\,x\in C\}$ exist and satisfy $0\leq\phi_1-\phi_2\leq 1$ and $\phi_1-1\leq\phi(\shl)\leq\phi_2$ for any line bundle $\shl$ on $C$. Up to a $\grp$-action, one can assume that $\phi_1=1$, then $0\leq\phi(\shl)\leq\phi_2\leq 1$.

Since $\sigma$ is non-numerical, one can choose a class $\ell$ in $\Pic^{\circ}(C)\subset\Pic(C)$ such that $Z(\ell)$ is contained in $\mathbb{U}\cup\br_{>0}$. Let $g>0$ be the genus of $C$, then the line bundle $\sho_C(gx+n\ell)$ admits a global section for any $x\in C$ and $n\in\bz$ due to the Riemann--Roch formula. It follows a non-trivial morphism $\shl\rightarrow\shl(gx+n\ell)$ and therefore $0\leq\phi(\shl)\leq\phi(\shl(gx+n\ell))$.

Suppose that $Z(\ell)\in \mathbb{U}$, then $Z(\shl(gx+n\ell))\notin\mathbb{U}\cup\br$ for $n$ small enough which contradicts the fact that $0\leq\phi(\shl(gx+n\ell))\leq 1$. 

Suppose that $Z(\ell)>0$, then the phase of any line bundle is $0$ as $\phi(\shl(gx+n\ell))$ can be arbitrary close to zero for $n\gg0$. On the other hand, one can use a similar argument and 
$$\phi(\shl(-gx-n\ell))\leq\phi(\shl)\leq 1$$
to show that any line bundle has phase $1$, a contradiction!
\end{proof}

In particular, combining with Theorem \ref{3-macri-main-result}, one sees immediately that

\begin{corollary}
One has $\Stab(\DCoh(C))=\Stab(C)\cong\bc\times\bh$.
\end{corollary}

The same result has been claimed as a special case of \cite[Theorem 4.6]{MP14} according to the remark after it.\footnote{Thanks Emanuele Macrì for informing this reference to the author.} However, the author would like to point out that, without Proposition \ref{3-main-1-1}, the last paragraph in \cite[p. 381]{MP14} is not necessary true for any locally-finite stability conditions.

\subsection{Classification of non-locally-finite stability conditions}
Here we show the second main result of this note about non-locally-finite stability conditions.

\begin{proposition}\label{3-main-1-2}
Up to a unique $\bc$-action, a non-locally-finite stability condition on $\DCoh(C)$ is equal to $\sigma_{\beta}$ in Example \ref{2_slope-stability-degenerate} for some $\beta\in\br-\bq$.
\end{proposition}

\begin{proof}
Choose a non-locally-finite numerical stability condition $\sigma$, one can assume that $\phi(\sho_x)=1$ and $Z(\sho_x)=-1$ for any $x\in C$ up to a unique $\bc$-action. Then there exists a line bundle with phase $0$ or $1$, otherwise one can argue as in \cite[Theorem 2.7]{Mac07} to see $\sigma\in\Stab(C)$. Also, for any $\sigma$-semistable coherent sheaf $\shg$ one has $0\leq\phi(\shg)\leq 1$ according to Corollary \ref{3-lemma-macri}.

Therefore, the central charge has the form
$$Z(r,d)=-d + \beta r$$
for some real number $\beta$, so $\schp((0,1])=\schp(1)\neq\cate{Coh}(C)$. One can see that a line bundle $\shl$ lies in $\schp(1)$ if and only if the usual slope $\mu(\shl)>\beta$ and is an object in $\schp(0)$ if and only if the usual slope $\mu(\shl)<\beta$. Then we \emph{claim}
$$\schp(0)\cap\cate{Coh}(C)=\schf_{\beta}\quad\textup{and}\quad \schp(1)\cap\cate{Coh}(C)=\scht_{\beta}$$
so that $\schp(1)=\schp((0,1])=\scha_{\beta}$ by Lemma \ref{2_fact-t-structure}.

We show the claim by doing induction on the rank. The rank $\leq1$ cases are true and we assume that the claim holds for coherent sheaves with rank $\leq n$. 

At first, one sees that the claim is true for non-slope-semistable locally free sheaves by taking the Harder--Narasimhan filtrations for slope stability. Next, any slope semistable locally free sheaf $\she$ on $C$ of rank $n+1$ must be $\sigma$-semistable.

Otherwise the Harder--Narasimhan filtration of $\she$ with respect to the stability condition $\sigma$, together with Lemma \ref{3-key-technical-lemma}, ensures a short exact sequence $0\rightarrow\shh\rightarrow\she\rightarrow\shg\rightarrow0$ of coherent sheaves such that $\shh\notin\schp(0)$ and $\shg$ is $\sigma$-semistable with $\phi(\shg)=0$. Here one uses the fact that $\phi(\shg)$ is contained in $[0,1]$ and the inductive assumption. Since $\she$ is slope semistable, one has $\mu(\shh)\leq\mu(\she)\leq\mu(\shg)<\beta$ and any slope semistable Harder--Narasimhan factor of $\shh$ has slope less than $\beta$. It follows that $\shh\in\schp(0)$ and a contradiction. Therefore, the claim is true for locally free sheaves of rank $n+1$ .

Consequently, one sees that a coherent sheaf $\shf$ of rank $n+1$ whose torsion-part $\sht$ is non-trivial cannot be in $\schp(0)$ by the non-trivial morphism $\sht\rightarrow\shf$ and $\phi(\sht)=1$. Since the claim holds for locally free sheaves with rank $n+1$, the sheaf $\shf$ is in $\schp(1)$ if and only if its torsion-free part satisfies $\mu_->\beta$. So we have finished the induction step and can conclude.

The claim also implies that $\beta$ is irrational as for each $r>0$ and $d\in\bz$ there exists a slope semistable sheaf with rank $r$ and degree $d$. Hence, one concludes.
\end{proof}

\begin{remark}
One should notice that the $\grp$-orbit of $\sigma_{\beta}$ is in fact equal to the $\bc$-orbit thereof. Indeed, an element $a+b\bm{i}$ in $\grp/\bc\cong\bh$ corresponds to 
$$M=\begin{pmatrix}
b&a\\
0&b^{-1}
\end{pmatrix}$$
whose action on stability conditions sends $\sigma_{\beta}$ to the pair $(\schp_{\beta},b^{-1}\cdot Z_{\beta}(r,d))$.
\end{remark}

\begin{remark}\label{3_P1-non-locally-finite}
On $\bp^1$, there are no non-locally-finite stability conditions by \cite{Oka06}, and the pair $\sigma_{\beta}=(\schp_{\beta},Z_{\beta})$ constructed in Example \ref{2_slope-stability-degenerate} is indeed a stability condition on $\DCoh(\bp^1)$ for every $\beta\notin\bz$ corresponding to the one described in \cite[Proposition 3.3 (3)]{Oka06} up to tensor with a line bundle and a $\bc$-action. In the complex manifold $\Stab(\bp^1)\cong\bc\times\bc$, these stability conditions are in the boundary of the geometric chamber $\Geo(\bp^1)\cong\bc\times\bh$. One also compares to \cite[Section 3.2]{Li17}.
\end{remark}

\begin{remark}
One should be able to establish some parallel results for triangulated categories analogous to $\DCoh(C)$. One can, for example, apply the argument of Proposition \ref{3-main-1-2} to classify non-locally-finite stability conditions on the bounded derived category $\DCoh(C_{FF})$ of an algebraic Fargues--Fontaine curve $C_{FF}$ based on the expositions in \cite[Section 2.3.4 and Section 2.4]{DJL24}.
\end{remark}

\subsection{The distances towards non-locally-finite stability conditions}
Even though there are no nice generalized metric on the set of non-locally-finite stability conditions, one has a generalized metric the set of all slicings \cite[Section 6]{Bri07}.

\begin{definition}\label{3_topology-slicing}
Let $\scht$ be a triangulated category, one defines a topology by
$$d(\schp_1,\schp_2)=\inf\{\epsilon\geq0\,|\,\schp_2(\phi)\subset\schp_1([\phi-\varepsilon,\phi+\varepsilon])\,\textup{ for all } \phi\in\br\}$$
on the set $\Slic(\scht)$ of all slicings on $\scht$.
\end{definition}

Since all the stability conditions on $\DCoh(C)$ are numerical, the set of all stability conditions can be endowed with the subspace topology from the product space
$$\Slic(\DCoh(C))\times\Hom(N(C),\bc)$$
according to \cite[Section 6]{Bri07}. In this case, one can make the following observations.

\begin{proposition}
Consider a smooth curve $C$ of positive genus and the following two cases:
\begin{itemize}
	\item[(1)] $\sigma_1=(\schp_1,Z_1)$ and $\sigma_2=(\schp_2,Z_2)$ are two non-locally-finite stability conditions on $\DCoh(C)$ that are not in the same $\grp$-orbit;
	\item[(2)] $\sigma_1=(\schp_1,Z_1)$ is a locally-finite stability condition on $\DCoh(C)$ and $\sigma_2=(\schp_2,Z_2)$ is a non-locally-finite stability condition on $\DCoh(C)$;
\end{itemize}
then one has $d(\schp_1,\schp_2)\geq \frac{1}{2}$ for either case. 

\begin{proof}
Here we only prove the first case, the second one is similar. Choose two non-locally-finite stability conditions $\sigma_1$ and $\sigma_2$. Then, up to a $\grp$-action, one can assume that $\sigma_1=\sigma_{\beta_1}$ and $\sigma_2=\lambda\cdot\sigma_{\beta_2}$ for some $\lambda\in\bc$ and $\beta_1\neq\beta_2$. Suppose in addition that $\beta_1<\beta_2$ without losing generality. Then one can find an semistable bundle $\she$ in $\schp_{\beta_1}(1)\cap\schp_{\beta_2}(0)$. Since $\sigma_2=\lambda\cdot\sigma_{\beta_2}$, one can write $\schp_2(\phi_0)=\schp_{\beta_2}(0)$ for some $\phi_0\in\br$ and then $\she\in\schp_{\beta_1}(1)\cap\schp_2(\phi_0)$. On the other hand, one can always find an object $A\in\schp_{\beta_1}(1)\cap\schp_{\beta_2}(1)=\schp_{\beta_1}(1)\cap\schp_2(1+\phi_0)$. It follows
$$d(\schp_1,\schp_2)\geq \max\{|1-\phi_0|,|\phi_0|\}\geq\frac{1}{2}$$
where all the equality holds if and only if $\phi_0=1/2$.
\end{proof}
\end{proposition}

It means that the subspace of non-locally-finite stability conditions is topologically a disjoint union of countably many $\bc$. These copy of $\bc$ do not attach to the complex manifold $\Stab(C)$. 

On the other hand, the distance $d(Z_{\beta},Z_{\alpha,\beta})$ in $\Hom(N(C),\bc)$ approaches zero when the index $\alpha$ approaches zero. So the non-locally-finite stability conditions can still be seen as certain kinds of boundary points for $\Stab(C)$ once some information on the slicing side is forgot.

\section{Boundary of the stability manifolds for smooth curves}

\subsection{Non-locally-finite stability conditions as boundary points}\label{4_Section4.1}
There are several ways to describe a boundary for the stability manifold in literature. Here we will see that non-locally-finite stability conditions on $\DCoh(C)$ are natural boundary points for a positive genus smooth curve $C$ from three different perspectives.

\subsubsection{Partial boundary of the naive closure}\label{4_Section4.1.1}
Image of the local homeomorphism
$$\pi\colon\Stab(C)\rightarrow\Hom_{\bz}(K(C),\bc),\quad (\schp,Z)\mapsto Z$$
in Theorem \ref{2_main-Bridgeland} is the open submanifold $\GL^+(2,\br)\subset\Hom_{\bz}(N(C),\bc)=\bc^2$. The complex structure on $\GL^+(2,\br)$ must be left invariant i.e.\ the left multiplication of each element in $\GL^+(2,\br)$ is holomorphic. Then one has an isomorphism 
$$\GL^+(2,\br)\cong \bc^\times\times\bh$$ 
as complex manifolds according to \cite[Theorem 2]{Sa81}. In fact, the left invariant complex structure on $\GL^+(2,\br)$ is unique up to a real number and one can choose the following representative
$$\GL^+(2,\br)\cong \bc^\times\times\bh,\quad 
\begin{pmatrix}
x_1&x_2\\x_3&x_4
\end{pmatrix}
\mapsto \left(\frac{1}{x_4+x_3\bm{i}},\frac{x_1-x_2\bm{i}}{x_3-x_4\bm{i}}\right)$$
of the unique class. Then the closure of $\GL^+(2,\br)$ is $\bar{\bc}^\times\times\bar{\bh}$ such that $\partial\bc^\times=\{\infty_1\}$ and $\partial\bh=\br\cup\{\infty_2\}$. So one has $\partial\GL^+(2,\br)=(\bar{\bc}^\times\cup\partial\bh)\cup(\{\infty_1\}\cup\bh)$.

The central charge $Z_{\beta}(r,d)=-d+\beta r$ of the stability condition $\sigma_{\beta}$ corresponds to the point $(\bm{i},-\beta)$ and the orbit $Z_{\beta}.\GL^+(2,\br)$ corresponds to $\bc^\times\times\{-\beta\}$. 

\begin{remark}
The above procedure can be seen as a special case of Bolognese \cite{Bol23}.
\end{remark}

\subsubsection{Partial boundary of the Thurston compactification}
Motivated by Thurston compactification of the Teichmüller space $\mathfrak{T}(\Sigma_g)$ of marked Riemann surfaces of genus $g>0$, a partial compactification for the quotient $\Stab(\scht)/\bc$ by the $\bc$-action is proposed in \cite{BDL23}. On smooth curves, it is studied by \cite{KKO24}.

At first, one defines the \emph{projective space} $\bp_{\geq0}^{\schs}$ as the quotient topological space
$$\bp_{\geq0}^{\schs}:=(\br_{\geq0}^{\schs}-\{0\})/\br_{>0}$$
for a subset $\schs$ of the set of isomorphism classes of objects in $\scht$. 

Then, for a suitable choice of $\schs$, one has a continuous map
$$\bp m\colon \Stab(\scht)/\bc\rightarrow\bp_{\geq0}^{\schs},\quad \sigma\mapsto[m_{\sigma}(E)]_{E\in\schs}$$
where $m_{\sigma}(E)$ is the sum of $|Z_{\sigma}(A_i)|$ for all Harder--Narasimhan factors $A_i$ of $E$. Suppose that the map $\bp m$ is a homeomorphism onto its image, one will obtain a compactification for $\Stab(\scht)/\bc$.

\begin{theorem}[{{\cite[Theorem 1.1]{KKO24}}}]
Consider a smooth curve $C$ of positive genus and choose $\schs=\{\sho_x,\sho_C,\sho_C(-y)\}$ for two points $x,y$ on $C$, then the following continuous map is homeomorphic onto the image
$$\bp m\colon\Stab(C)/\bc\rightarrow\bp_{\geq0}^{\schs},\quad \sigma\mapsto[m_{\sigma}(\sho_x):m_{\sigma}(\sho_C):m_{\sigma}(\sho_C(-y))]$$
and its closure is homeomorphic to the closed hyperbolic disk. 
\end{theorem}

More precisely, one has $\Stab(C)/\bc\cong\bh$ and the point $\beta+\alpha\bm{i}\in\bh$ corresponds to the stability condition $\sigma_{\alpha,\beta}$ with central charge $Z_{\alpha,\beta}(r,d)=-d+(\beta+\alpha\bm{i})r$. Then
$$\bp m(\sigma_{\alpha,\beta})=[1:\sqrt{\beta^2+\alpha^2}:\sqrt{(\beta+1)^2+\alpha^2}]$$
and a boundary point $[1:\beta:\beta+1]$ with $\beta$ irrational can be seen as image of the non-locally-finite stability condition $\sigma_{\beta}$ in Theorem \ref{3-main-1-2} under $\bp m$. It gives the subset $\br-\bq$ of $\partial\bh=\br\cup\{\infty\}$.

\begin{remark}
The continuous map $\bp m$ is not injective for $C=\bp^1$ and the same $\schs$ but the image of $\bp m(\Stab(\bp^1)/\bc)$ gives a partial compactification of $\Geo(\bp^1)/\bc$ according to \cite{KKO24}. This should be compared with Remark \ref{3_P1-non-locally-finite}.
\end{remark}

\subsubsection{Boundary points of the global dimension closure}
In \cite{IQ23,Qiu23}, Ikeda and Qiu introduce a notion of dimension for a given stability condition. This notion is useful in investigate certain properties of the stability manifold, see \cite{FLM23,FLLQ23}. Here we relate it to the boundary of $\Stab(C)$.

\begin{definition}
The \defi{global dimension} of a pre-stability condition $\sigma=(\schp,Z)$ is 
$$\gldim(\sigma):=\sup\{\phi_2-\phi_1\,|\,\Hom_{\scht}(A_1,A_2)\neq0\textup{ for }A_i\in \schp(\phi_i)\}$$
which ranges from $0$ to $\infty$.
\end{definition}

The global dimension of a stability condition is preserved by the $\bc$-action and the $\Aut(\scht)$-action by definition. Moreover, one has 

\begin{theorem}
Let $C$ be a smooth curve of genus $g$, then
\begin{itemize}
	\item If $g=1$, then $\gldim(\sigma)=1$ for any stability condition $\sigma\in\Stab(C)$;
	\item If $g\geq 2$, then $2>\gldim(\sigma)>1$ for any stability condition $\sigma\in\Stab(C)$;
\end{itemize}
In fact, one has $\sup\{\gldim(\sigma)\,|\,\sigma\in\Stab(C)\}=2$ for smooth curves of genus $g\geq 2$.

\begin{proof}
Thanks to \cite[Theoreom 5.16]{KOT21}, it remains to show that $\gldim(\sigma)<2$ for any $\sigma\in\Stab(C)$ and the supremum of global dimension is $2$. Up to a $\bc$-action, one can assume that $\sigma=\sigma_{\alpha,\beta}$ and conclude that $\gldim(\sigma_{\alpha,\beta})<2$ as the homological dimension of $\cate{Coh}(C)$ is $1$. The supremum is due to $\gldim(\sigma_{\alpha,0})\geq 1+\phi(\omega_C)-\phi(\sho_C)\rightarrow 2$ for $\alpha>0$ small enough
\end{proof}
\end{theorem}

The supremum cannot be reached for curves of genus $g\geq 2$. Since 
$$\gldim\colon \Stab(C)\rightarrow(1,2)$$
is continuous \cite[Lemma 5.7]{IQ23}, it is natural to expect the global dimension of a point in the boundary of $\Stab(C)$ is $1$ or $2$. It can be verified for non-locally finite stability conditions.

\begin{proposition}
A non-locally-finite stability condition on $\DCoh(C)$ for a smooth curve $C$ has global dimension $1$ once the genus of $C$ is $1$, and has global dimension $2$ once the genus of $C$ at least $2$.

\begin{proof}
It suffices to compute $\gldim(\sigma_{\beta})$ by Proposition \ref{3-main-1-2}. Due to \cite[Remark 5.5]{IQ23}, it is equal to the global dimension of the abelian category $\scha_{\beta}$. One only needs to check $\Hom(X,Y[2])$ by construction. Once $g(C)=1$, an easy computation shows
$$\Hom(X,Y[2])=\Hom(\shh^{-1}(X),\shh^0(Y)[1])=\Hom(\shh^0(Y),\shh^{-1}(X))^\vee=0$$
so one has $\gldim(\sigma_{\beta})=1$. 

Once $g(C)\geq 2$, one can find some line bundle $\shl$ on $C$ such that $\mu(\shl)<\beta$ and $\mu(\omega_C\otimes\shl)>\beta$. So one has $\shl[1],\omega_C\otimes\shl\in\scha_{\beta}$ but $\Hom(\shl[1],\omega_C\otimes\shl[2])=\bc$.
\end{proof}
\end{proposition}

\subsection{Other boundary points and CLSY weak stability}
The non-locally-finite stability conditions do not make the whole boundary in many cases. To described other boundary points, some weaker versions of stability conditions are introduced such as \cite{BPPW22,CLSY25}. David informs that their definition in \cite{BPPW22} will change \cite{Plo25}, so only the notion in article \cite{CLSY25} will be discussed here.

\begin{definition}
A \defi{CLSY weak stability condition} $(\schp,Z,\{\phi_A\})$ on a triangulated category $\scht$ consists of a pre-stability condition $\sigma=(\schp,Z)$ and a collection of real numbers $\{\phi_A\}_{A\in S_{\sigma}}$ indexed by
$$S_{\sigma}=\{A\in\schp((0,1])\,|\,A\neq0\textup{ and }Z(A)=0\}$$
such that $\phi(A)=\phi_A$ for any semistable object $A\in S_{\sigma}$ and for any short exact sequence $0\rightarrow K_1\rightarrow K\rightarrow K_2\rightarrow0$ in $\schp((0,1])$, one of the following two inequalities
\begin{align*}
	\phi_{K_1}\geq\phi_K\geq\phi_{K_2}\quad\textup{ and }\quad\phi_{K_1}\leq\phi_K\leq\phi_{K_2}
\end{align*}
holds, where $\phi_A$ is defined by $Z(A)=m(A)e^{\bm{i}\pi\phi_A}$ for  $A\in\schp((0,1])-S_{\sigma}$.
\end{definition}

The two group actions on the set of all pre-stability conditions on $\scht$ preserve the subset of all CLSY weak stability conditions on $\scht$.

\begin{example}
A stability condition is a CLSY weak stability condition.
\end{example}

\begin{example}\label{4_example-CLSY-weak-slope}
The slope stability on a smooth surface $S$ determines a numerical CLSY weak stability condition with heart $\cate{Coh}(S)$, such that $\phi(\sht)=1$ for any coherent sheaf $\sht$ with zero-dimensional support. This CLSY weak stability condition is not a stability condition as the central charge of $\sht$ is zero.
\end{example}

The pre-stability condition $\sigma_{\beta}$ in Example \ref{2_slope-stability-degenerate} and $\sigma'_{\beta}$ in Example \ref{2_slope-stability-degenerate-alternative} can be seen as the natural boundary points of $\Stab(C)$ as in Section \ref{4_Section4.1}.

\begin{example}\label{4_example-CLSY-weak-boundary}
Let $\beta\in\br\cup\{\infty\}$, then the pre-stability condition $\sigma_{\beta}=(\schp_{\beta},Z_{\beta})$ becomes a CLSY weak stability condition by claiming that $\phi(A)=1$ for any $A$ with $A\in\ker(Z)\cap\scha_{\beta}$. The same construction works for the pre-stability condition $\sigma'_{\beta}$ for any $\beta\in\br\cup\{\infty\}$.
\end{example}

However, this definition is too loose such that for each $\beta\in\bq\cup\{\infty\}$ one can find a lot of extra CLSY weak stability conditions.

\begin{example}\label{4_example-CLSY-weak-boundary-abundant}
Let $\beta\in\bq\cup\{\infty\}$ and $C$ be a smooth curve of positive genus, then one can define a CLSY weak stability condition $\sigma_{\beta,t}=(\schp_{\beta,t},Z_{\beta},\{\phi(A)\})$ for each real number $t\in[0,1]$ such that for $t\neq1$ one has
\begin{align*}
\schp_{\beta,t}(1)&=\scha_{\beta}^{\circ}:=\{A\in\scha_{\beta}\,|\,Z_{\beta}(A)\neq0\}\\
\schp_{\beta,t}(t)&=\{A\in\scha_{\beta}\,|\,Z_{\beta}(A)=0\}
\end{align*}
and for $t=1$ one defines $\sigma_{\beta,1}:=\sigma_{\beta}$. One notices that $\sigma'_{\beta}=\sigma_{\beta,0}$.
\end{example}

In general, it is hard to control the CLSY weak stability condition with degenerated central charge and one can define infinitely many of them on a slicing. 

\begin{example}
Let $\sigma=(\schp,Z)$ be a stability condition, then for any non-zero object $A$ in $\schp((0,1])$ one has $\phi_{\sigma}(A)\in(0,1]$ determined by $Z(A)=m(A)e^{\bm{i}\pi\phi_{\sigma}(A)}$. Then one can define a new CLSY weak stability condition $(\schp_f,O,\{\phi_f(A)\})$ for the zero map $O\colon K(\scht)\rightarrow\bc$ and a monotonic function $f\colon(0,1]\rightarrow(0,1]$ such that $\schp_f(\phi):=\schp(f(\phi))$ and $\phi_f(A):=f(\phi_{\sigma}(A))$.
\end{example}

\subsection{The weak stability conditions}
The previous subsection indicates one needs a stronger notion of weak stability conditions to shape the boundary. Here we suggest a definition for positive genus smooth curves.

\begin{definition}
A \defi{weak stability condition} $\sigma$ on a triangulated category $\scht$ is a CLSY weak stability condition $(\schp,Z,\{\phi_A\})$ such that $Z$ is non-trivial and for any short exact sequence 
$$0\rightarrow K_1\rightarrow K\rightarrow K_2\rightarrow0$$
in $\schp((0,1])$, then one of the following inequalities 
\begin{align*}
\phi_{K_1}>\phi_K>\phi_{K_2},\quad\phi_{K_1}<\phi_K<\phi_{K_2},\quad\phi_{K_1}=\phi_K=\phi_{K_2}
\end{align*}
holds, where $\phi_A$ is defined as for CLSY weak stability conditions.
\end{definition}

The two group actions on the set of all pre-stability conditions on a given triangulated category $\scht$ preserve the subset of all weak stability conditions on $\scht$.

\begin{example}
The CLSY weak stability condition in Example \ref{4_example-CLSY-weak-slope} is not a weak stability condition by the short exact sequence $0\rightarrow\shi_x\rightarrow\sho_S\rightarrow\sho_x\rightarrow0$.
\end{example}

\begin{example}
A stability condition is a weak stability condition. The CLSY weak stability conditions in Example \ref{4_example-CLSY-weak-boundary} are weak stability conditions.
\end{example}

Moreover, we want more constrains on phases to make sure that the phase function will not go wild when one takes limit from $\Stab(C)$ to the boundary.

\begin{definition}
A weak stability condition $\sigma=(\schp,Z,\{\phi(A)\})$ on the bounded derived category $\DCoh(X)$ of a smooth variety $X$ is called \defi{regular} if $\schp((0,1])$ does not contain objects with trivial numerical class.
\end{definition}

\begin{example}\label{4_example-weak-stability-condition}
A numerical stability condition on $\cate{D}^b(X)$ is a regular weak stability condition. The weak stability conditions in Example \ref{4_example-CLSY-weak-boundary} are regular weak stability conditions.
\end{example}

The two group actions on the set of all pre-stability conditions on $\DCoh(X)$ preserve the subset of all regular weak stability conditions on $\DCoh(X)$.

\begin{proposition}
Up to a unique $\bc$-action, any regular weak stability condition on $\DCoh(C)$ is equal to a weak stability condition in Example \ref{4_example-weak-stability-condition} for a smooth curve $C$ of positive genus.

\begin{proof}
Let $\sigma$ be a weak stability condition on $\DCoh(C)$. According to Proposition \ref{3-main-1-1}, its central charge $Z$ factors through the numerical Grothendieck group $N(C)$. Suppose that the image of $Z$ has rank $2$, then one can apply the argument of Theorem \ref{3-macri-main-result} to see that $\sigma$ is indeed a locally finite stability condition. Suppose that the image of $Z$ has rank $1$, then one has two cases.

Assume that $Z(0,1)=0$, then after a unique $\bc$-action, one can assume $Z(r,d)=r$ and that the phase of all line bundles on $C$ is $0$. Then, for a given $x\in C$ one has either $\phi(\sho_x)=0$ or $\phi(\sho_x)=1$ according to Corollary \ref{3-lemma-macri} and the definition. Here one notices that
$$0\rightarrow \sho_x\rightarrow\sho_C(-x)[1]\rightarrow\sho_C[1]\rightarrow0$$
is a short exact sequence in $\schp((0,1])$ when $\sho_x\in\schp((0,1])$. Since $\sigma$ is regular, the phase of any skyscraper sheaves on $C$ must be the same. In this case, one reduces to two cases: either $\phi(\sho_x)=1$ for any $x\in C$ and then $\sigma=\sigma_{\infty}$, or $\phi(\sho_x)=0$ for any $x\in C$ and then $\sigma=\sigma'_{\infty}$. 

Assume that $Z(0,1)\neq0$, then after a unique $\bc$-action, one can assume that $Z(r,d)=-d+\beta r$ for some $\beta\in\br$ and that the phase of all skyscraper sheaves on $C$ is $1$. Let $r_0$ be the minimal rank for a semistable locally free sheaf with slope $\beta$, then using Corollary \ref{3-lemma-macri} one can argue inductively as in Proposition \ref{3-main-1-2} to show
\begin{align*}
\schp(0)\cap\cate{Coh}_{<r_0}(C)&=\schf_{\beta}\cap\cate{Coh}_{<r_0}(C)=\schf'_{\beta}\cap\cate{Coh}_{<r_0}(C)\\
\schp(1)\cap\cate{Coh}_{<r_0}(C)&=\scht_{\beta}\cap\cate{Coh}_{<r_0}(C)=\scht'_{\beta}\cap\cate{Coh}_{<r_0}(C)
\end{align*}
where $\cate{Coh}_{< r_0}(C)$ denotes the category of coherent sheaves with rank $< r_0$. Then one can show that a semistable sheaf $\she_{r_0}$ with slope $\beta$ and rank $r_0$ is $\sigma$-semistable as in Proposition \ref{3-main-1-2} again when $r_0>1$ or use Corollary \ref{3-macri-main-result} when $r_0=1$. Then one has either $\phi(\she_{r_0})=0$ or $\phi(\she_{r_0})=1$ as there always exists a short exact sequence of coherent sheaves
$$0\rightarrow\shh\rightarrow\she_{r_0}\rightarrow\shg\rightarrow0$$
such that $\shh\in\schf_{\beta}$ and $\shg\in\scht'_{\beta}$. It means that either
$$\schp(0)\cap\cate{Coh}_{\leq r_0}(C)=\schf_{\beta}\cap\cate{Coh}_{\leq r_0}(C)\quad\textup{and}\quad \schp(1)\cap\cate{Coh}_{\leq r_0}(C)=\scht_{\beta}\cap\cate{Coh}_{\leq r_0}(C)$$
or
$$\schp(0)\cap\cate{Coh}_{\leq r_0}(C)=\schf'_{\beta}\cap\cate{Coh}_{\leq r_0}(C)\quad\textup{and}\quad \schp(1)\cap\cate{Coh}_{\leq r_0}(C)=\scht'_{\beta}\cap\cate{Coh}_{\leq r_0}(C)$$
depending on $\phi(\she_{r_0})$. Since $\sigma$ is regular, any $\sigma$-semistable coherent sheaf with slope $\beta$ has same phase with the sheaf $\she_{r_0}$ and one can argue inductively as above and the method in Proposition \ref{3-main-1-2} to show that $\sigma=\sigma_{\beta}$ once one has $\phi(\she)=1$ and $\sigma=\sigma'_{\beta}$ once one has $\phi(\she)=0$. The details are left to the readers.
\end{proof}
\end{proposition}

In general, one can get infinitely many weak stability conditions on a positive genus smooth curve $C$ for each $\beta\in\bq\cup\{\infty\}$ and any assignment of phase $0$ or $1$ to stable coherent sheaves with slope $\beta$ on $C$. The associated bounded $t$-structures can be described analogous to \cite[Proposition 5.1]{BK06}.

\section{Stability conditions on singular curves}
In this section, we will investigate the stability conditions on the bounded derived category of a singular curve. 

\subsection{The Grothendieck group of singular curves}
At first, one needs to understand the Grothendieck group $K(C)$. One has a homomorphism
$$\ch\colon K(C)\rightarrow \CH_{\bullet}(C)\otimes_{\bz}\bq$$
according to the Hirzebruch--Riemann--Roch theorem (see, for example, \cite{BFM75}), whose base change
$$\ch_{\bq}\colon K(C)\otimes_{\bz}\bq\rightarrow \CH_{\bullet}(C)\otimes_{\bz}\bq$$
is a group isomorphism. In particular, the kernel $T(C)$ of the surjection
$$(\rank,\deg)\colon K(C)\rightarrow \bz^{\oplus 2}:=N(C)$$ 
contains finite rational spans of zero-cycles on $C$ with degree zero. So up to multiple with an integer, any class in $T(C)$ becomes a degree zero Weil divisor. Unlike the smooth case, one cannot always find a canonical line bundle for a Weil divisor. However, one always has a choice.

\begin{proposition}\label{5-techinical-lemma}
Consider a class $\ell$ in $T(C)$, then there exists an integer $n$ and a degree zero line bundle $\shl$ such that the class of $[\shl]-[\sho_C]$ in $K(C)$ equals to $n\ell$.

\begin{proof}
By the comments just before the proposition, one can take an integer $n>0$ such that $n\ell$ is the class of a Weil divisor on $C$ say $\sum n_i[x_i]$. It suffices to show that, for a point $x_i$, one can find a line bundle $\shl$ such that the class of $[\shl]-[\sho_C]$ in $K(C)$ equals to $[x_i]$. One takes the normalization $f\colon\tilde{C}\rightarrow C$ and choose a point $\tilde{x}_i$ over $x_i$, then $\tilde{x}_i$ is rational equivalence to some zero cycle $Z$ on $\tilde{C}$ whose support is disjoint from $f^{-1}(x_i)$. The pushforward $f_*\tilde{x}_i$ is rational equivalence to $f_*Z$ so by definition $x_i=f_*Z+H$ where $H$ is the Cartier divisor determined by a rational function $h$ on $C$. Since the support of $f_*Z$ is disjoint from $x_i$, one has $x_i=H$ on an open neighborhood $U$ of the point $x_i$. Then one can define a Cartier divisor $D$ on $C$ by $D|_U=H$ and $D_{C-\{x_i\}}=1$. It provides the desired line bundle $\shl$.
\end{proof}
\end{proposition}

This argument is attributed to Qing Liu based on \cite[Page 599]{ColT05}. A general statement can be found in \cite[Theorem 4.5 and 6.5]{GLL13}.

\subsection{A connected component of the stability manifold}
Similar to the smooth case, one can consider the following definition:

\begin{definition}
One defines the space $\Geo^{\dag}(C)$ to be the $\grp$-orbit of the locally finite stability condition in Example \ref{2_slope-stability}.
\end{definition}

 As before, it is diffeomorphic to $\bc\times\bh$ and is a complete metric space by \cite{Woolf12}. So it is closed in every metric space containing it. In particular, $\Geo^{\dag}(C)\subset\Stab(\DCoh(C))$ is a closed submanifold.

According to \cite[Proposition 8.1]{Bri07}, the complex linear topology on $\Stab(\DCoh(C))$ is locally given by the topology in Definition \ref{3_topology-slicing}

Now we are prepared to prove Theorem \ref{1_main-3}. The basic idea is to show that the stability conditions in $\Geo^{\dag}(C)$ cannot deform away from $\Geo^{\dag}(C)$.

\begin{proposition}
The subspace $\Geo^{\dag}(C)$ is a connected component of $\Stab(\DCoh(C))$.

\begin{proof}
It suffices to show that for any $\sigma^{\dag}\in\Geo^{\dag}(C)$ one has
$$B_{\epsilon}(\sigma^{\dag})=\{\sigma\in\Stab(\DCoh(C))\,|\,d(\sigma,\sigma^{\dag})<\epsilon\}\subset \Geo^{\dag}(C)$$
for some small $\epsilon>0$. Up to a unique $\grp$-action, it suffices to show it for the stability condition $\sigma_0=(\schp_0,Z_0)$ constructed in Example \ref{2_slope-stability}.

Choose a point $\sigma=(\schp,Z)$ in $B_{\epsilon}(\sigma_0)$, at first we claim that the central charge $Z$ factors through $N(C)$ via the rank and degree functions. Otherwise, one can choose a class $\ell$ in the kernel $T(C)$ such that $Z(-\ell)\in \mathbb{U}\cup\br_{>0}$. Up to multiplication with a positive integer, one can assume that $\sho(\ell)$ is a degree zero line bundle according to Proposition \ref{5-techinical-lemma}. In this case, $\sho_C(n\ell)$ is a line bundle on $C$ any $n\in\bz$ and is stable with respect to $\sigma_0$. Then one has
$$\sho_C(n\ell)\in\schp[\frac{1}{2}-\epsilon,\frac{1}{2}+\epsilon]$$
for any $n\in\bz$ as $\phi_{\sigma}(\sho_C(n\ell))=1/2$ and $d(\sigma,\sigma_0)<\epsilon$. It is impossible for sufficiently small $\epsilon$. Therefore, one concludes $Z$ factors through $N(C)$.

In particular, combining with the fact that $\Geo^{\dag}(C)$ is a connected space, the subspace $\Geo^{\dag}(C)$ is contained in a connected component $\Stab^{\dag}(\DCoh(C))$ of the complex manifold $\Stab(\DCoh(C))$. Then one has a subspace $V^{\dag}\subset\Hom(K(C),\bc)=\bc^2$ and a local homeomorphism
$$\pi\colon\Stab^{\dag}(\DCoh(C))\rightarrow V^{\dag},\quad (\schp,Z)\mapsto Z$$
by Theorem \ref{2_main-Bridgeland}. The local homeomorphism restricts to the universal covering
$$\pi\colon\Geo^{\dag}(C)\cong\grp\rightarrow \GLp(2,\br)$$
so once $\sigma=(\schp,Z)$ lies in $B_{\epsilon}(\sigma_0)-\Geo^{\dag}(C)$, it satisfies $Z\notin\GLp(2,\br)$. Since the subspace $\GL^+(2,\br)\cong\bc^\times\times\bh$ is connected by Section \ref{4_Section4.1.1} and, thanks to \cite[Lemma 6.4]{Bri07}, the local homeomorphism $\pi$ is bijective on $B_{\epsilon}(\sigma_0)$ for $\epsilon$ small enough, one reduces to show that the central charge of any $\sigma=(\schp,Z)$ in $ B_{\epsilon}(\sigma_0)$ cannot belong to the boundary of $\GLp(2,\br)$.

It is true because any boundary point corresponds to a linear function $Z\colon N(C)\rightarrow\bc$ with image isomorphic to $\br$. Then one can always find some $\phi\in[0,1]$ such that $\schp[\phi-\varepsilon,\phi+\varepsilon]$ is empty but in the meanwhile, according to \cite{Bal99}, one can find a semistable bundle on $C$ with phase sufficiently closed to $\phi$. It contradicts the fact that $d(\sigma,\sigma_0)<\epsilon$. So we are done.
\end{proof}
\end{proposition}

Similar to the positive genus smooth case, the pair $\sigma_{\beta}=(\schp_{\beta},Z_{\beta})$ is a non-locally-finite stability condition on $\DCoh(C)$ for any $\beta\in\br-\bq$ and can be seen as natural boundary points for the component $\Geo^{\dag}(C)$. However, Proposition \ref{3-main-1-2} does not apply directly in this case.

\begin{remark}
The connected component $\Geo^{\dag}(C)$ is always preserved by the action of the group $\Aut(\DCoh(C))$ on $\Stab(\DCoh(C))$. In fact, suppose that $C$ is strict Cohen--Macauley or is Gorenstein with ample of anti-ample canonical bundle, then one concludes by the identification
$$\Aut(\DCoh(C))\cong\Aut(C)\rtimes(\Pic(C)\times\bz)$$
in \cite[Lemma 2.7 and The Proof of Theorem 1.1]{Sp23} and \cite[Proposition 6.18]{Ball11}. 

Otherwise $C$ is Gorenstein with neither ample or anti-ample canonical bundle. Then the arithmetic genus of $C$ is $1$ and one has $\Stab(\DCoh(C))=\Geo^{\dag}(C)$ by \cite{BK06}. 
\end{remark}

\subsection{Geometric stability conditions on singular curves}
This subsection is devoted to the second part of Theorem \ref{1_main-3}: the geometric stability conditions are contained in $\Geo^{\dag}(C)$. The critical point is that a (projective and integral) curve $C$ admits at most Cohen--Macauley singularities. 

\begin{definition}
Let $X$ be a variety, then a complex of coherent sheaves is called \defi{perfect} if it is quasi-isomorphic to a bounded complex of locally free sheaves.
\end{definition}

Since the curve $C$ is Cohen--Macauley, its dualizing complex $\omega_C$ concentrates at degree zero and one has the following duality result for perfect complexes \cite{Hart}.

\begin{theorem}[Serre duality]
One has a functorial isomorphism
$$\Hom(A,B[n])\cong\Hom(B,A\otimes\omega_C[1-n])^\vee$$
for objects $A,B$ in $\DCoh(C)$ such that one of them is perfect.
\end{theorem}

It is the fundamental tool for the first step. Here one also recalls that $\scht_{\infty}\subset\DCoh(C)$ is the subcategory of torsion sheaves on $C$.

\begin{proposition}\label{5_heart-of-geometric-stable}
Consider a geometric stability condition $\sigma=(\schp,Z)$ on $C$ such that all skyscraper sheaves have phase $1$, then $\schp(1)=\scht_{\infty}$ and $\schp((0,1])=\cate{Coh}(C)$.

\begin{proof}
Choose a stable object $A$ in $\schp(1)$, then $A\cong\sho_x$ for some $x\in C$. Otherwise one has
$$\Hom(\sho_x,A[n])=0\quad\textup{and}\quad \Hom(A,\sho_x[n])=0$$
for any $n\leq0$ and $x\in C$. The object $\sho_x$ is perfect when $x$ is regular, so one has
$$\Hom(A,\sho_x[n])=\Hom(A,\sho_x\otimes\omega_C[n])=\Hom(\sho_x,A[1-n])=0$$
by the Serre duality, for any $n\geq1$ and regular point $x\in C$. Here one also uses the fact that the canonical sheaf $\omega_C$ is locally free at the Gorenstein points. So by \cite[Lemma 5.3]{BM02} one concludes that $A$ is supported on singular points. Moreover, the condition
$$\Hom(A,\sho_x[n])=0,\quad\forall n\leq0,x\in C$$
ensures that $\shh^q(A)=0$ for $q\geq0$ according to \cite[Proposition 5.4]{BM02}. Choose the minimal $m$ such that $\shh^m(A)\neq0$, then $m\leq-1$ and the spectral sequence 
$$E_2^{p,q}=\Hom(\sho_x,\shh^q(A)[p])\Rightarrow\Hom(\sho_x,A[p+q])$$
and $\Hom(\sho_x,A[m])=0$ for any $x\in C$ implies that $\Hom(\sho_x,\shh^m(A))=0$ for any $x\in C$. It means that $\shh^m(A)$ is torsion-free, contradicting that the support of $\shh^m(A)$ is discrete.

Choose an object $E$ in $\schp((0,1))$, then one has
$$\Hom(\sho_x,E[n])=0\quad\textup{and}\quad \Hom(E,\sho_x[n-1])=0$$
for any $n\leq0$ and $x\in C$. One sees as before that the restriction of $E$ on the regular locus of $C$ concentrates to a locally free sheaf at degree zero. Moreover, one has $\shh^q(E)=0$ for $q\geq1$ and the support of $\shh^q(E)$ is contained in the singular locus for $q\leq -1$. Then one argues as before to see that $\shh^q(E)=0$ for $q\leq -1$ and the coherent sheaf $\shh^0(E)$ is torsion-free.

In conclusion, one has $\schp(1)=\scht_{\infty}$ and $\schp((0,1))=\schf_{\infty}$ so that one can conclude the desired identification $\schp((0,1])=\cate{Coh}(C)$.
\end{proof}
\end{proposition}

\begin{proposition}
Consider a geometric stability condition $\sigma=(\schp,Z)$ on $\DCoh(C)$ such that all the skyscraper sheaves have phase $1$ and $\schp((0,1])=\cate{Coh}(C)$, then $\sigma$ has to be numerical.

\begin{proof}
Since $\schp((0,1])=\cate{Coh}(C)$, all the line bundles on $C$ are semistable. Otherwise, one can find a Harder--Narasimhan filtration of a line bundle $\shl$ in $\schp((0,1])$ and in particular a non-trivial injection of sheaves $\shg\hookrightarrow\shl$ whose cokernel is not torsion, which is absurd.

Suppose otherwise that $\sigma$ is not numerical, then according to Proposition \ref{5-techinical-lemma} one can find a class $\ell$ in the kernel $T(C)=\ker(K(C)\rightarrow N(C))$ such that one can find a line bundle $\shm$ on $C$ with class $[\sho_C]+\ell$ in $K(C)$ and $Z(\ell)>0$. Then one concludes by Proposition \ref{3-main-1-1}.
\end{proof}
\end{proposition}

\begin{proposition}
Consider a numerical geometric stability condition $\sigma=(\schp,Z)$ on $\DCoh(C)$ such that all the skyscraper sheaves have phase $1$ and $\schp((0,1])=\cate{Coh}(C)$, then $\sigma$ is in $\Geo^{\dag}(C)$.

\begin{proof}
It suffices to show that the central charge $Z$ of $\sigma$ is contained in $\GLp(2,\br)$. At first, one can assume that $Z(r,d)=-d+(\beta+\alpha\bm{i}) r$ up to a constant. Since $Z(r,d)$ for any $d\in\bz$ and $r\geq0$ is contained in $\mathbb{U}$, one has $\alpha>0$ as well.
\end{proof}
\end{proposition}

The second part of Theorem \ref{1_main-3} follows, since after a $\bc$-action one can assume that all the skyscraper sheaves have phase $1$ for a geometric stability condition.

\subsection{Numerical stability conditions on non-rational curves}
It is proved in \cite{FLZ22} that any numerical stability condition on a smooth projective variety with finite Albanese morphism is in fact geometric. The key tool to control the behaviors of stability conditions is the following statement. Here one notices that, according to \cite[Corollary 9.17]{LO10}, every element in $\Aut(\DCoh(X))$ is a \emph{kernel} autoequivalence in the sense of \cite{Pol07} for an irreducible projective scheme $X$.

\begin{lemma}[{{\cite[Theorem 3.5.1]{Pol07}}}]
Consider an irreducible projective scheme $X$ and a bounded $t$-structure with noetherian heart on $\DCoh(X)$. Assume that a connected smooth algebraic group $G$ acts on $\DCoh(X)$ by autoequivalences, then the bounded $t$-structure is invariant under this action.
\end{lemma}

It is known that the connected algebraic group $\Pic^0(C)$ is smooth (see \cite{BLR}), so any stability condition in $\Stab(C)$ is preserved by the group action by $\Pic^0(C)$. Here one should notice that the heart of a general element in $\Stab(C)$ is noetherian according to \cite[Corollary 5.0.2]{AP06}.

\begin{proposition}
Consider a non-rational singular curve $C$, then the skyscraper sheaves are all stable with respect to every numerical stability condition on $\DCoh(C)$.

\begin{proof}
Fix an element $\sigma$ in $\Stab(C)$ and pick any closed point $x\in C$. Let $A_1,\dots,A_n$ be the Jordan--Hölder factors in the Harder--Narasimhan filtration for $\sho_x$ with respect to $\sigma$. Then the lemma implies that $A_i\otimes\shl\cong A_i$ for any $\shl\in\Pic^0(C)$. Moreover, for any $q$ one has
$$\pi^*\shh^q(A_i)\otimes\pi^*\shl\cong \pi^* \shh^q(A_i)$$
where $\pi\colon \tilde{C}\rightarrow C$ is the normalization. Since the homomorphism $\pi^*\colon\Pic^0(C)\rightarrow\Pic^0(\tilde{C})$ is surjective, the support of $\pi^*\shh^q(A_i)$ is finite by the argument of \cite[Theorem 2.13]{FLZ22}. It follows that the support of $\shh^q(A_i)$ is also finite, as the normalization $\pi$ is a finite morphism.

So one can conclude by \cite[Lemma 2.14]{FLZ22}. In fact, its argument works whenever one can find a stable factor of $\sho_x$ supported at $x$.  
\end{proof}
\end{proposition}

The next step is to prove an analogy of \cite[Proposition 2.9]{FLZ22} for singular curves. It could be used to generalize \cite[Proposition 2.9]{FLZ22} for some singular varieties.

\begin{proposition}
Consider a singular curve $C$ and a numerical stability condition $\sigma$ on $\DCoh(C)$ such that every skyscraper sheaf on $C$ is stable with respect to $\sigma$, then the stability condition $\sigma$ is geometric.

\begin{proof}
One notices that \cite[Lemma 2.10]{FLZ22} works for a general projective variety, so one can apply the argument of \cite[Proposition 2.9]{FLZ22} to see that all the skyscraper sheaves at smooth points on $C$ have the same phase. Then, up to a $\grp$-action, one can assume that the phase of those skyscraper sheaves is $1$. In this case, the stable objects in $\schp(1)$ are either skyscraper sheaves or shifts $\sho_{x_i}[2n_i]$ for some singular points $x_1,\dots,x_k$ and non-zero integers $n_1\leq n_2\leq\dots\leq n_k$ by a suitable modification of the arguments for Proposition \ref{5_heart-of-geometric-stable}. 

Now choose $E\in\schp((0,1))$ and denote $D=\{x_1,\dots,x_k\}$. Then one has
$$\Hom(\sho_x,E[n])=0\quad\textup{and}\quad \Hom(E,\sho_x[n-1])=0$$
$$\Hom(\sho_x,E[n-2n_i])=0\quad\textup{and}\quad \Hom(E,\sho_x[n-1+2n_i])=0$$
for any $n\leq0$, any point $x\in C-D$, and any $i$. Similar to Proposition \ref{5_heart-of-geometric-stable}, one can see that the restriction of $A$ to $C-D$ concentrates in a torsion-free sheaf at degree zero and the support of $\shh^i(E)$ is contained in $D$ for $i\neq0$. It follows that $\shh^0(E)$ is supported on the whole of $C$.

Suppose that $n_1>0$, then one has $\shh^q(E)=0$ for $q\geq 1$ due to \cite[Proposition 5.4]{BM02}. Then the minimal integer $m$ such that $\shh^{-m}(E)\neq0$ is zero, so that the spectral sequence
$$E_2^{p,q}=\Hom(\shh^{-q}(E),\sho_x[p])\Rightarrow\Hom(E,\sho_x[p+q])$$
ensures that $\Hom(E,\sho_x)=\Hom(\shh^0(E),\sho_x)\neq0$ for any $x\in C$. It contradicts $n_k>0$. 

Suppose that $n_1<0<n_k$, then one can restrict to the open subset $U\subset C$ which is the complement of the discrete set $D_0=\{x_i\in D\,|n_i<0\}$. Since $U_0$ is quasi-projective, the argument for the last case still works: one can obtain $\shh^q(E_U)=0$ for $q\geq 1$ and
$$\Hom(E_U,\sho_x)=\Hom(\shh^0(E_U),\sho_x)\neq0$$
for any $x\in U$, which is impossible as $n_k>0$.

Suppose that $n_k<0$, then one has $\shh^q(E)=0$ for $q\leq -1$ and $\shh^0(E)$ is torsion-free using the argument for Proposition \ref{5_heart-of-geometric-stable}. Also, the spectral sequence
$$E_2^{p,q}=\Hom(\sho_x,\shh^q(E),[p])\Rightarrow\Hom(\sho_x,E[p+q])$$
implies that $\Hom(\sho_x,\shh^0(E)[1])=0$ for any $x\in D$. It can never happen. Indeed, one always has a short exact sequence $0\rightarrow\shh^0(E)\rightarrow\pi_*\she\rightarrow\sht\rightarrow0$ where $\pi\colon \tilde{C}\rightarrow C$ is the resolution of singularities, $\she$ is the quotient of $\pi^*\shh^0(E)$ by its torsions, and $\sht$ is a torsion on $C$. Once $x$ is not in the support of $\sht$, one has $$\Hom(\sho_x,\shh^0(E)[1])=\Hom(\sho_x,\pi_*\she[1])=\Hom(L\pi^*\sho_x,\she[1])\neq0$$ 
using the fact that $x\times_C\tilde{C}$ is finite. Otherwise, one concludes by $\Hom(\sho_x,\sht)\neq0$.

In conclusion, we show that neither case is possible and the phase of every skyscraper sheaf has to be $1$ i.e. $\sigma$ is geometric.
\end{proof}
\end{proposition}

Combining the two propositions, one obtains the desired statement.

\begin{theorem}
Consider a non-rational singular curve $C$, then a numerical stability condition on $\DCoh(C)$ is geometric. In particular, one has $\Stab(C)=\Geo^{\dag}(C)$.
\end{theorem}

In this case, one can describe a partial boundary of $\Stab(C)$ similar to the case for smooth curves discussed in Section 4. However, it is still difficult to classify the non-locally-finite numerical stability conditions on singular curves due to the lack of a parallel version of Corollary \ref{3-lemma-macri}.

Also, it is not clear that whether an arbitrary stability condition on a singular curve is numerical. However, one should notice that it is true for any rational singular curve $C$ since in this case one has $K(C)=\bz^{\oplus 2}$ using \cite[Lemma 3.1]{BK06}.

\section*{Acknowledgments}
The author would like to thank his advisors Laura Pertusi and Paolo Stellari for their constant supports, thank Anand Deopurkar and David Ploog for exchanges on compactifications at an early stage, thank Tom Bridgeland and Emanuele Macrì for helpful suggestions, and thank Zhiyu Liu for beneficial conversations. Also, the author would like to thank the anonymous referees for comments.

\end{document}